
\documentstyle{amsppt}
\magnification=1200

\hoffset=-0.5pc
\vsize=57.2truepc
\hsize=38truepc
\nologo
\spaceskip=.5em plus.25em minus.20em

 \define\anderone{1}
 \define\atibottw{2}
 
 \define\atiysemi{4}
 
\define\bottsega{6}
 \define\bouskan{7}

\define\donalkro{10}

\define\dwykanon{12}

  \define\modus{16}
\define\modustwo{17}
    \define\srni{18}

\define\huebjeff{20}

\define\jeffrtwo{22}
\define\jeffrthr{23}
  
  \define\kantwo{25}

\define\karshone{28}

\define\maclbotw{31}

\define\puppeone{36}
\define\shulmone{37}

\define\weinstwe{39}

\define\wittesix{42}

 \define\atiyahfo{43}
\define\atiboton{44}
\define\bergever{45}
  \define\bottone{46}
 \define\bottsame{47}
\define\botshust{48}
    \define\chen{49}
 \define\chentwo{50}
\define\doldpupp{51}
\define\guhujewe{52}
\define\huebkan{53}

\define\akingone{54}
\define\newstone{55}
\define\sourithr{56}
\define\sourifiv{57}

\define\witteele{58}

\noindent
dg-ga/9506006
\bigskip\noindent
\topmatter
\title
Extended moduli spaces
and the Kan construction. II.
\\
Lattice gauge theory
\endtitle
\author Johannes Huebschmann{\dag}
\endauthor
\affil
Max Planck Institut f\"ur Mathematik
\\
Gottfried Claren Str. 26
\\
D-53 225 BONN
\\
huebschm\@mpim-bonn.mpg.de
\endaffil
\thanks
{{\dag} The author carried out this work in the framework of the
VBAC research group of EUROPROJ.}
\endthanks
\abstract{Let $Y$ be a CW-complex with a single 0-cell, let $K$ be its Kan
group, a free simplicial group whose realization is a model for the space
$\Omega Y$ of based loops on $Y$, and let $G$ be a compact, connected Lie
group. We carry out an explicit purely finite dimensional construction of
generators of the equivariant cohomology of the geometric realization of the
cosimplicial manifold $\roman{Hom}(K,G)$ and hence, in view of earlier results,
of the space $\roman{Map}^o(Y,BG)$ of based maps from $Y$ to the classifying
space $BG$ of $G$ where $G$ acts on $BG$ by conjugation. For a smooth manifold
$Y$, this may be viewed as a rigorous approach to lattice gauge theory, and we
show that it then yields, (i) when {$\roman{dim}(Y)=2$,} equivariant de Rham
representatives of generators of the equivariant cohomology of twisted
representation spaces of the fundamental group of a closed surface including
generators for moduli spaces of semi stable holomorphic vector bundles on
complex curves so that, in particular, the known structure of a stratified
symplectic space results; (ii) when {$\roman{dim}(Y)=3$,} equivariant
cohomology generators including the Chern-Simons function; (iii) when
{$\roman{dim}(Y) = 4$,} the generators of the relevant equivariant cohomology
from which for example Donaldson polynomials are obtained by evaluation against
suitable fundamental classes corresponding to moduli spaces of ASD
connections.}
\endabstract
\date{June 13, 1995} \enddate
\endkeywords
\endsubjclass

\endtopmatter
\document
\leftheadtext{Johannes Huebschmann}
\rightheadtext{Lattice gauge theory}

\beginsection Introduction

In this paper we pursue further
the purely finite dimensional approach to gauge theory
began in \cite\huebkan:
Let $Y$ be a CW-complex with a single 0-cell,
let $K$ be its Kan group
\cite\kantwo,
a free simplicial group
whose realization is a model for the space
$\Omega Y$ of based loops on $Y$, and
let $G$ be a compact and connected Lie group.
In this paper we carry out an explicit construction
of
the generators of the $G$-equivariant de Rham cohomology
of
the realization $|\roman{Hom}(K,G)|$
 of the cosimplicial manifold
$\roman{Hom}(K,G)$
and hence,
in view of the main result in our paper \cite\huebkan,
of
the space
$\roman{Map}^o(Y,BG)$
of based maps from $Y$ to the classifying space
$BG$ of $G$
where $G$ acts on $BG$ by conjugation;
we thereby exploit the fact that the chains on the simplicial
nerve of $K$ yield a model for
the chains of $Y$.
Since
$\roman{Hom}(K,G)$
is a smooth finite dimensional
manifold in each cosimplicial degree,
cf. \cite\huebkan,
the construction is
{\it purely finite dimensional:
every de Rham form
will be constructed on a piece of finite dimension\/}.
A concise statement is given in Theorem 7.1 below.
The finite dimensional pieces
belong to the cosimplicial manifold
$\roman{Hom}(K,G)$.
Integration
then carries these forms to forms
on the
realization $|\roman{Hom}(K,G)|$.
This requires  a suitable interpretation of forms
on mapping spaces.
Using the theory of {\it differentiable space\/}
\cite\chen, \cite\chentwo\ or, what amounts to the same,
that of {\it \lq\lq diffeological\rq\rq\  space\/}
(\lq\lq espace diff\'eologique\rq\rq)
\cite\sourithr, \cite\sourifiv,
forms on mapping spaces
admit a purely
finite dimensional interpretation
in terms of what are called
{\it plots\/}
\cite\chen, \cite\chentwo\
or
\lq\lq {\it plaques\/}\rq\rq\
\cite\sourithr, \cite\sourifiv\
and do {\it not\/} require infinite dimensional techniques.
For a smooth manifold $Y$, our construction
of forms
may be viewed as a rigorous approach to
lattice gauge theory,
whereby plots admit a natural interpretation
as (equivariant) {\it families of principal bundles with connection\/};
see Section 5 below for details.
\smallskip
We offer three applications; they
may be viewed as
classical  topological field theory
constructions:
At first we show that the construction yields, when
{$\roman{dim}(Y)=2$,}
explicit
equivariant de Rham representatives
of
equivariant generators of the cohomology of
moduli spaces of
twisted representation spaces
of the fundamental group of a closed surface;
in particular,
this yields
the structure of a stratified symplectic space
on such a moduli space
already obtained by other means
\cite\modus, \cite\srni, \cite\huebjeff.
We expect that part of what is said in \cite\witteele\
can be understood within our framework.
It is  worthwhile pointing out, though, that
even for the case of a bundle on a closed surface
$\Sigma$,
the present more general
construction involving a model for the full loop space
rather than merely a presentation of the fundamental group
of the surface \cite\modus, \cite\modustwo, \cite\huebjeff,
\cite\jeffrtwo, \cite\jeffrthr\
goes beyond
earlier constructions:
The realization
$|\Cal H|$ of
$\Cal H = \roman{Hom}(K\Sigma,G)$
contains the spaces
of based gauge equivalence classes
of {\it all\/}
central Yang-Mills connections
\cite\atibottw,
not just those which correspond
to the absolute minimum or, equivalently,
to projective representations
of the fundamental group $\pi$ of $\Sigma$,
and hence the space
$|\Cal H|$
comes with a kind of Harder-Narasimhan
filtration,
cf. Section 2 of \cite\huebkan.
The latter cannot
be obtained from the earlier extended moduli space
constructions.
Perhaps information about the multiplicative
structure of the cohomology of moduli spaces
can be derived from the  models we shall construct
below or from variants thereof.
Secondly,
when {$\roman{dim}(Y)=3$,}
we obtain equivariant  cohomology
generators
including
an explicit expression for the Chern-Simons function
on our model of the space of based
gauge equivalence classes of connections,
thereby answering a question
raised by {\smc Atiyah}
in \cite\atiysemi\
where he comments
on a possible combinatorial approach to
the path integral quantization
of the Chern-Simons function.
Thirdly, when {$\roman{dim}(Y) = 4$,}
we obtain
the generators of the equivariant cohomology of the appropriate space
from which
for example
Donaldson polynomials are obtained
by evaluation against suitable fundamental classes
corresponding to moduli spaces of ASD connections.

\smallskip
Our construction
is rigid in the sense that it gets away with
various choices made in the earlier approaches;
the theory, admittedly technically a bit complicated,
will take care of itself,
{\it no\/}
choices of appropriate data
must be made except that of various chains
representing
certain homology classes,
and the occurrence of the homotopy
operator on forms
in
the cited references
will get its natural explanation
in terms of
a realization procedure involving
integration of forms;
see Section 5 below for details.
Our approach is vastly more general
than those in
\cite\modus, \cite\modustwo,
\cite\huebjeff,
\cite\jeffrtwo,
\cite\jeffrthr\
since
it applies
to a bundle over
an arbitrary smooth compact manifold
via a cell decomposition or triangulation
as explained above.
Formally it is not even necessary to know that
the simplicial group
we are working with
arises from a smooth manifold;
we shall therefore expose the theory
for an arbitrary simplicial group or groupoid.
In this way we arrive at a kind of gauge theory
over arbitrary
CW-complexes.
By means of the
simplicial groupoid
constructed in
\cite\dwykanon\
for an arbitrary connected
simplicial set
the present approach
can be extended
to arbitrary
connected simplicial complexes,
in particular,
to triangulated smooth manifolds.
\smallskip
Our models for the space of
gauge equivalence classes of connections
involve
classical low dimensional topology notions
such as {\it identity among relations\/}
(Section 3 of \cite\huebkan)
and
{\it universal quadratic group\/}
(Section 4 of \cite\huebkan);
this somewhat establishes a link between classical algebraic
topology and the more recent gauge theory
developments in low dimensions.
We expect that our models will also prove
useful for various calculations
recently done in quantum cohomology
and that related
finite dimensional constructions
may be applied to other gauge theory
situations.
\smallskip
Some historical comments
about the origin of the present purely finite dimensional techniques
may be in order:
Extending an approach by
{\smc Karshon}
\cite\karshone,
A. Weinstein \cite\weinstwe\ constructed
a closed equivariant 2-form on
(the smooth part)
of certain spaces of homomorphisms
$\roman{Hom}(\pi,G)$
from the fundamental
group $\pi$ of a closed surface
to a Lie group $G$
with a biinvariant metric
and showed by techniques from equivariant
cohomology \cite\atiboton\
that this 2-form descends to
(the non-singular part of)
$\roman{Rep}(\pi,G)$.
In
\cite\modus, \cite\huebjeff, and \cite\jeffrtwo,
Weinstein's method has been refined
so as to
yield a smooth finite dimensional
symplectic manifold with a hamiltonian action
so that the space of representations
and more general twisted
versions thereof arise
by symplectic reduction;
this approach has been extended thereafter
in \cite\modustwo\ and
\cite\guhujewe\
to more general planar groups than just
surface groups so that for example
moduli spaces of parabolic bundles
can be successfully treated.
Another generalization in
\cite\jeffrthr\
yields
explicit representatives for
{\smc Newstead's}
generators
\cite\newstone\
for
various
moduli spaces
over a surface;
in the algebro-geometric context,
these arise as moduli spaces of
certain semi stable holomorphic vector bundles
on complex curves.
The present paper
gives such a construction
for
an arbitrary
gauge theory situation.
It may be viewed as the
\lq\lq grand unified theory\rq\rq\
searched for by A.\,Weinstein
in \cite\weinstwe.
Our {\it principal innovation\/}
is to replace
the bar construction
of a discrete group
coming into play in \cite\karshone,
\cite\weinstwe,
and in the subsequent papers
\cite\modus, \cite\modustwo,
\cite\huebjeff,
\cite\jeffrtwo,
\cite\jeffrthr,
\cite\guhujewe,
by the simplicial nerve of the Kan group $K$
on $Y$
so that
we can handle
the space of based gauge equivalence classes
of connections on an  arbitrary
principal bundle
with compact structure group
on a {\it general\/} manifold $Y$.
\smallskip
Any unexplained notation is the same
as
that in our papers \cite\modus\
and \cite\huebkan.
Details about cosimplicial spaces
may be found in
 \cite\bottsega\
and \cite\bouskan.
Topological spaces will be assumed endowed
with the compactly generated topology throughout.
\smallskip
I am indebted
to S. Bauer, H. J. Baues, and D. Puppe for discussions,
and to J. Stasheff for a number of most helpful
comments on a draft of the manuscript.
The paper has been written during a stay
at the Max Planck Institut at Bonn.
I wish to express my
gratitude to it
and to its director Professor
F. Hirzebruch for hospitality and support.

\medskip\noindent{\bf 2. Forms on spaces of representations}
\smallskip\noindent
Let $\Pi$ be a finitely generated groupoid, for example a group, and
write $(C_\natural(\Pi),\partial_\natural)$
and $(C^\natural(\Pi),\delta^\natural)$ for the
complexes of
normalized chains
and cochains, respectively,
on its {\it nerve\/}
$N\Pi$ or {\it inhomogeneous reduced
normalized
bar construction}.
We use the
dummy symbol $\natural$
to distinguish
bar resolution and hence group or groupoid (co)-homology degree
from form degree which will be written $*$.
Further, let $G$ be a connected Lie group;
the extension of the construction to be given below
to general non-connected Lie groups will
be studied elsewhere.
View $G$ as a groupoid with a single object,
and consider the space
$\roman H = \roman {Hom}(\Pi,G)$.
This space is not necessarily  smooth
at every point, and the interpretation of de Rham forms will
in general require some care.
However in the present paper we shall only need the special case
where $\Pi$ is {\it free\/}
so that $\roman H$ amounts to a product of finitely many
copies of $G$.
\smallskip
Equivariant de Rham forms
on
$\roman H$
may be constructed in the following way:
Given a $k$-tuple
$[x_1|x_2|\dots|x_k]$ of elements of
$\Pi$, $k \geq1$,
and an equivariant de Rham form
$\alpha \in \Omega^{i,j}_G(G^k)$, $ i, j \geq 0$, the evaluation map
$$
E_{[x_1|x_2|\dots|x_k]}
\colon
\roman{Hom}(\Pi,G)
@>>>
G^k,
\quad
\phi \mapsto (\phi(x_1),\dots,\phi(x_k)),
\tag2.1
$$
yields the form
$$
E^*_{[x_1|x_2|\dots|x_k]}(\alpha)
\in
\Omega^{i,j}_G(\roman{Hom}(\Pi,G)).
\tag2.2
$$
This construction can be formalized in the following way:
\smallskip
Let $k\geq 0$,
and consider
the
differential graded algebra
$$
\Omega^*(\roman{Hom}(\Pi,G) \times \Pi^k)
=
\Omega^{*}(\roman{Hom}(\Pi,G))
\otimes C^{k}(\Pi).
\tag2.3
$$
The evaluation map
$
E
$
from
$\roman{Hom}(\Pi,G) \times \Pi^k$
to
$G^k$
is compatible with the
obvious $G$-actions
and
induces
a morphism
$$
E^{*,*} \colon
(\Omega_G^{*,*}(G^k),d,\delta_G) @>>>
(\Omega_G^{*,*}(\roman{Hom}(\Pi,G));d,\delta_G)
 \otimes  C^{k}(\Pi)
\tag2.4
$$
of equivariant de Rham algebras.
Moreover,
as $k$ varies,
these maps assemble to a morphism
$$
\left(\Omega_G^{*,*}(G^{\natural}); d,\delta_G,\delta^{\natural}\right)
@>>>
\left(\Omega_G^{*,*}(\roman{Hom}(\Pi,G)); d,\delta_G\right)
\otimes (C^{\natural}(\Pi),\delta^{\natural})
\tag2.5
$$
of tricomplexes;
in a given tridegree $(i,j,k)$, it
goes from
$\Omega_G^{i,j}(G^k)$
to
\linebreak
$\Omega_G^{i,j}(\roman{Hom}(\Pi,G))
\otimes C^{k}(\Pi)$.
For each bar complex degree $k$,
pairing with chains in
$C_{k}(\Pi)$, we obtain the graded bilinear pairing
$$
\langle\cdot,\cdot\rangle
\colon
\left((\Omega_G^{*,*}(\roman H); d,\delta_G)\otimes
C^{k}(\Pi)\right)
\otimes C_{k}(\Pi)
@>>>
\left(\Omega_G^{*,*}(\roman H); d,\delta_G\right)
\tag2.6
$$
which is compatible with the operators
$d$ and $\delta_G$
and,
for every $u \in \Omega_G^{*,*}(\roman H)\otimes C^{k}(\Pi)$
and every
$v \in
C_{k+1}(\Pi)$,
satisfies
$$
\langle u, \partial_\natural v \rangle
=
(-1)^{k+1}\langle \delta^{\natural} u, v \rangle,
\quad
\tag2.7
$$
where
the right-hand side refers to
(2.6) for $k+1$ rather than $k$;
here the
sign
$(-1)^{k+1}$
is forced by the
Eilenberg-Koszul convention
for the differential on a Hom-complex.
Combining (2.6) with (2.5)
and abusing the notation
$\langle\cdot,\cdot\rangle$
slightly,
we then obtain the pairing
$$
\langle\cdot,\cdot\rangle
\colon
\left(\Omega_G^{*,*}(G^{\natural}); d,\delta_G,\delta^{\natural}\right)
\otimes (C_\natural(\Pi),\partial_\natural)
@>>>
\left(\Omega_G^{*,*}(\roman H); d,\delta_G\right)
\tag2.8
$$
which is compatible with the operators
$d$ and $\delta_G$
and, moreover, satisfies
$$
\langle Q, \partial_\natural c \rangle
=
(-1)^{k+1}\langle \delta^{\natural} Q, c \rangle,
\quad
Q \in \Omega_G^{*,*}(G^k), \
c \in
C_{k+1}(\Pi),
\tag2.9
$$
whatever $k \geq 0$.
Thus pairing a form
$Q$
in
$\Omega_G^{*,*}(G^{\natural})$
against a chain $c$ in
$C_\natural(\Pi)$,
we obtain
the form
$\langle Q,c\rangle$ in
$\Omega_G^{*,*}(\roman H)$.
We shall need an explicit expression
for the value
$D\langle Q,c\rangle$
in terms of $Q$ and $c$
of the total differential
$D$
on
the right-hand side of (2.8).
There is no real obstacle
to calculating
this value
in terms of
the pairing (2.8)
and the operators
$d,\delta_G,
\delta^{\sharp}$, and $\partial^{\sharp}$,
but since (2.8)
does not behave as a pairing
of chain complexes
for the operators
$\delta^{\sharp}$ and $\partial^{\sharp}$,
cf. (2.9), this calculation
is somewhat of a mess.
The cure is provided by
an extension of
the construction
which leads to the formula (2.15) below:
Recall that, for an arbitrary differential graded coalgebra $C$
with diagonal $\Delta$
and arbitrary ground ring $R$, --- in fact, we could take an arbitrary
differential graded algebra here ---
the cap pairing
$
\cap
$ from
$\roman{Hom}(C,R)
\otimes C$
to
$C$
is given by the composite
$$
\roman{Hom}(C,R)
\otimes C
@>{\roman {Id} \otimes \Delta}>>
\roman{Hom}(C,R)
\otimes C
\otimes C
@>{\roman{ev} \otimes \roman{Id}_C}>>
R \otimes C
@>{\cong}>>
C
$$
where \lq\lq ev\rq\rq\ denotes the evaluation pairing.
When we take for $C$ the inhomogeneous reduced  normalized
bar construction of $\Pi$, we obtain the cap pairing
from
$(C^\natural(\Pi),\delta^\natural)
\otimes
(C_\natural(\Pi),\partial_\natural)
$
to
$(C_\natural(\Pi),\partial_\natural)
$
inducing on homology the cap pairing
$
\cap
$
from
$\roman H^{\natural}(\Pi)
\otimes
\roman H_{\natural + \ell}(\Pi)
$
to
$\roman H_{\ell}(\Pi)$, for $\ell \geq 0$.
Tensoring the identity morphism with
the cap pairing
yields an extension
$$
\aligned
\roman{Id} \otimes \cap
&\colon
(\Omega_G^{*,*}(\roman H); d,\delta_G)\otimes
(C^{\natural}(\Pi),\delta^{\natural})
\otimes (C_{\natural}(\Pi),\partial_\natural)
\\
&\quad
@>>>
\left(\Omega_G^{*,*}(\roman H); d,\delta_G\right)
\otimes (C_{\natural}(\Pi),\partial_\natural)
\endaligned
\tag2.10
$$
of (2.6) above
which is compatible with all the operators
coming into play and hence
induces
a pairing
$$
\left|(\Omega_G^{*,*}(\roman H); d,\delta_G)\right|
\otimes
(C^{\natural}(\Pi),\delta^{\natural})
\otimes
(C_{\natural}(\Pi),\partial_\natural)
@>>>
|(\Omega_G^{*,*}(\roman H); d,\delta_G)|
\otimes (C_{\natural}(\Pi),\partial_\natural)
\tag2.11
$$
of the chain complexes
resulting from totalization,
which we write $| \cdot |$.
Recall that the total differential on
$\left|(\Omega_G^{*,*}(\roman H); d,\delta_G)\right|$
is simply the sum $d + \delta_G$.
Finally, when we combine (2.10) with (2.5),
we obtain the pairing
$$
\langle\cdot,\cdot\rangle
\colon
\left(\Omega_G^{*,*}(G^{\natural}); d,\delta_G,\delta^{\natural}\right)
\otimes (C_\natural(\Pi),\partial_\natural)
@>>>
\left(\Omega_G^{*,*}(\roman H); d,\delta_G\right)
\otimes (C_\natural(\Pi),\partial_\natural)
\tag2.12
$$
which is compatible with all the operators
and induces
a pairing
$$
\langle\cdot,\cdot\rangle
\colon
\left|(\Omega_G^{*,*}(G^{\natural}); d,\delta_G,\delta^{\natural})
\right|
\otimes
(C_{\natural}(\Pi),\partial_\natural)
@>>>
\left|
(\Omega_G^{*,*}(\roman H); d,\delta_G)
\right|
\otimes (C_{\natural}(\Pi),\partial_\natural)
\tag2.13
$$
of the chain complexes
resulting from totalization.
We remind the reader that,
for every $(i,j,k)$,
on the homogeneous
component
$\Omega_G^{i,j}(G^k)$,
the total differential
$d_G$ on
$\left|
(\Omega_G^{*,*}(G^{\natural}); d,\delta_G,\delta^{\natural})
\right|$
is given by
$$
d_G = d + \delta_G +(-1)^{i+j} \delta^{\natural}.
\tag2.14
$$
The compatibility property
of (2.13) means that, when
$D$ refers to the tensor product differential on
the right-hand side
$\left|(\Omega_G^{*,*}(\roman H); d,\delta_G)\right|
\otimes (C_{\natural}(\Pi),\partial_\natural)$
of (2.13),
for $Q
\in |\Omega_G^{*,*}(G^{\natural})|$
and
$c \in C_{\natural}(\Pi)$,
$$
D\langle Q,c\rangle
=
\langle d_GQ,c\rangle
+ (-1)^{|Q|}
\langle Q,\partial_\natural c\rangle
\tag2.15
$$
where $|Q|$ denotes the total degree of $Q$.
Notice in a given quadruple degree $(i,j,k,k+\ell)$, (2.12) goes from
$
\Omega_G^{i,j}(G^{k})
\otimes C_{k+\ell}(\Pi)
$
to
$\Omega_G^{i,j}(\roman H)
\otimes C_{\ell}(\Pi)$.
\smallskip
The pairing (2.12) and hence (2.13)
is natural, in fact {\it covariant\/}, in the variable $\Pi$
but notice that $\Pi$ also occurs in $\roman H = \roman{Hom}(\Pi,G)$
so that the forms $\Omega_G^{*,*}(\roman H)$
are also covariant in $\Pi$.
\smallskip
The total complex
$
\left|
(\Omega_G^{*,*}(G^{\natural}); d,\delta_G,\delta^{\natural})\right|
=
\left(|\Omega_G^{*,*}(G^{\natural})|,d_G\right)
$
inherits a structure of differential graded algebra
in the following way:
For each pairs $(i,j)$ and
$(i',j')$
of bidegrees
and for each $k,k'$,
consider the canonical pairing
$$
\Omega_G^{i,j}(G^k)
\otimes
\Omega_G^{i',j'}(G^{k'})
@>>>
\Omega_G^{i+i',j+j'}(G^k \times G^{k'});
$$
it amounts to the dual of the
{\it Alexander-Whitney\/} map
for the usual bar construction.
These pairings induce the searched for
structure of differential graded algebra.
It is natural in terms of the data.
\smallskip
The
differential graded algebra
$\left(|\Omega_G^{*,*}(G^{\natural})|,d_G\right)
$
computes the equivariant
real cohomology
algebra
of the classifying space $BG$ for $G$
where $G$ acts on $BG$ via conjugation.
To recall what this cohomology looks like,
we assume henceforth $G$ compact;
the general case may as usual be reduced
to this one by taking a maximal compact subgroup.
Let $Ig$ be the graded algebra
of invariant polynomials on $g$,
where $g$ is endowed with degree 2
as usual;
it is well known to be itself a finitely generated
polynomial algebra.
Inspection of the Serre spectral sequence
for the Borel construction
$EG \times_G BG$
shows at once that
the equivariant cohomology algebra
$\roman H_G^{*}(BG)$
of $BG$ is isomorphic to
$Ig
\otimes
Ig$.
In particular,
every
class in $\roman H^*(BG)$
has an
equivariantly closed representative in the total complex
$\left(|\Omega_G^{*,*}(G^{\natural})|,d_G\right)$,
that is,
the restriction mapping
from
$
\left(|\Omega_G^{*,*}(G^{\natural})|,d_G\right)$
to
$\left|(\Omega^{*}(G^{\natural}),d,\delta^\natural)\right|
$
induces a surjection
from
$
\roman H_G^*(BG)
$
to
$\roman H^*(BG)$
on cohomology.
\smallskip
Explicit generators arise as follows:
Take the realization of the nerve $NG$
of $G$ as a model for the classifying space $BG$,
and let $Q$ be an invariant degree $r$ polynomial on $g$.
{\smc Shulman's} simplicial Chern-Weil
construction
\cite\shulmone,
\cite\bottone,
\cite\botshust,
applied to the universal simplicial principal
$G$-bundle $(\xi_0, \xi_1,\dots )$
over the simplicial space
$NG$,
yields forms
$$
Q^{r,r} \in \Omega^r(G^r),
\ Q^{r+1,r-1} \in \Omega^{r+1}(G^{r-1}),
\ \dots \ ,
\ Q^{2r-1,1} \in
\Omega^{2r-1}(G),
\tag2.16
$$
and
the sum
$Q^{r,r} + \dots + Q^{2r-1,1}$
is a closed element of
$\left|(\Omega^{*}(G^{\natural});d,\delta^{\natural})\right|$
which represents the class $[Q] \in \roman H^{2r}(BG)
(= \roman H^{2r}(NG))$
arising from $Q$.
More precisely,
for each $q\geq 1$, the Maurer-Cartan forms
yield a connection on the corresponding (trivial) principal
$G$-bundle $\xi_q\colon G^{q+1} \times \Delta_q \to G^q \times \Delta_q$
having curvature
$F_q \in \Omega^2(G^q \times \Delta_q, \roman{ad}(\xi_q))$, and, for
$1 \leq q \leq r$,
$$
Q^{2r-q,q} = \int_{\Delta_q}Q(F_q)  \in \Omega^{2r-q}(G^q).
$$
Note that
$Q^{2r-1,1}$
is a closed form on $G$
representing the generator of
$\roman H^{2r-1}(G)$
which transgresses to $[Q]$.
\smallskip
As observed in \cite\jeffrthr,
the equivariant Chern-Weil construction
\cite\bergever\
yields explicit equivariant extensions of these forms:
Given
a Lie group $H$ and
an arbitrary $H$-equivariant
principal $G$-bundle
$\xi \colon P \to M$,
for an $H$-equivariant  connection on
$\xi$ with connection from $\vartheta \in \Omega^1(P,g)^H$,
define the {\it moment\/}
$\mu = \mu_\vartheta \in \Omega^{2,0}(M,\roman{ad}(\xi))$
of the connection by
$$
\mu \colon
h \to \Omega^0(M,\roman{ad}(\xi))
=
C^{\infty}(P,g)^G,
\quad
\mu(X) = \vartheta(X_P)
$$
where $X_P$ denotes the vector field on $P$ induced by $X$.
Then
an invariant degree $r$ polynomial $Q$  on $g$
determines the closed form
$$
Q(F+\mu)
= \widetilde Q^{0,2r} +  \widetilde Q^{2,2r-2}+ \dots + \widetilde Q^{2r,0}
\in  |\Omega_H^{*,*}(M)|^{2r}
$$
where $\widetilde Q^{i,j} \in \Omega_H^{i,j}(M)$.
When we apply this to the principal $G$-bundle
$\xi_q$ with $H=G$ acting by conjugation,
with the notation $\mu_q
\in \Omega^{2,0}(G^q\times \Delta_q,\roman{ad}(\xi_q))$
for the corresponding moment,
we obtain
the closed form
$$
Q(F_q+\mu_q)
= \widetilde Q^{0,2r} +  \widetilde Q^{2,2r-2}+ \dots + \widetilde Q^{2r,0}
\in  |\Omega_G^{*,*}(G^q\times \Delta_q)|^{2r}
$$
where $Q^{i,j} \in \Omega_G^{i,j}(G^q\times \Delta_q)$,
and integration
yields the forms
$$
\alignat 2
Q^{0,2r-q,q} &= \int_{\Delta_q}\widetilde Q^{0,2r}  \in \Omega^{0,2r-q}(G^q),
&&
\\
Q^{2,2r-2-q,q} &= \int_{\Delta_q}
\widetilde Q^{2,2r-2}  \in \Omega^{2,2r-2-q}(G^q),
&&
\\
&\cdots
\\
Q^{2r-q,0,q} &= \int_{\Delta_q}\widetilde Q^{2r-q,q} \in \Omega^{2r-q,0}(G^q),
&&
\text{\ if $q$ is even,}
\\
Q^{2r-q-1,1,q} &= \int_{\Delta_q}
\widetilde Q^{2r-q-1,q+1}  \in \Omega^{2r-q-1,1}(G^q),
\quad
&&
\text{\ if $q$ is odd.}
\endalignat
$$
For an invariant polynomial $Q$ on $g$ of degree $r$,
write
$$
\Omega_Q=\Sigma Q^{2i,j,q}, \quad 2i+j+q = 2r, \quad q \leq 2i+j, \quad q \leq
r;
$$
this is a closed element of
$\left|(\Omega^{*,*}(G^{\natural});d,\delta_G\delta^{\natural})\right|$
representing a class $[\Omega_Q] \in \roman H_G^{2r}(BG)$.

\proclaim{Theorem 2.17}
For every invariant polynomial $Q$ on $g$ of degree $r$,
the class $[\Omega_Q] \in \roman H_G^{2r}(BG)$
restricts to the class
$[Q] \in \roman H^{2r}(BG)$
arising from $Q$.
Furthermore,
when $Q$ runs through a set of polynomial
generators of $Ig$,
the classes
$[\Omega_Q]$
together with the elements $Q$ viewed
as elements of
$\Omega^{*,0}(G^0)$
constitute a set of polynomial generators
of
$
\roman H_G^*(BG)
=
\roman H^*\left|
(\Omega^{*,*}(G^{\natural});d,\delta_G,\delta^{\natural})\right|.
$
\endproclaim

\demo{Proof}
The first statement is immediate.
The
\lq\lq Furthermore\rq\rq\
clause is an immediate formal consequence thereof. \qed
\enddemo

For example, let $Q$  be an invariant symmetric bilinear
form $\cdot$ on $g$,
so that $r=2$.
The above
construction then yields
$$
\align
Q^{0,3,1} \in
\Omega_G^{0,3}(G),
\quad
Q^{2,1,1} \in
\Omega_G^{2,1}(G),
\quad &\text{\ for\ } q = 1,
\\
Q^{0,2,2} \in
\Omega_G^{0,2}(G\times G),
\quad
Q^{2,0,2} \in
\Omega_G^{2,0}(G\times G),
\quad &\text{\ for\ } q = 2,
\endalign
$$
and their sum
is
a
closed
4-form in the total complex
$\left|(\Omega^{*,*}(G^{\natural});d,\delta_G,\delta^\natural)\right|$.
Actually it may be shown that
the term
$Q^{2,0,2}$
is irrelevant and may be dropped.
The element
$Q^{0,3,1}$ is the fundamental 3-form
on $G$
constructed
by E. Cartan.
\smallskip
The
singular cochains $C^*(G)$ of $G$ constitute a Hopf algebra,
the requisite diagonal map being induced from
the
multiplication
mapping on $G$ by means of the shuffle map,
and it is well known and classical that
the cobar construction on $C^*(G)$
yields a model for the (singular) cochains on $BG$.
The bar de Rham bicomplex $(\Omega^*(G^\natural);d,\delta^\natural)$
serves as a replacement for the {\it cobar\/} construction
on the differential graded algebra
$\Omega^*(G)$
of forms on $G$
which is not available in the strict sense;
while the multiplication mapping of $G$ induces a map
from $\Omega^*(G)$ to
$\Omega^*(G\times G)$ we cannot
algebraically
project down the latter to
$\Omega^*(G) \otimes \Omega^*(G)$
in such a way that a coalgebra structure on
$\Omega^*(G)$ results.
The bar de Rham bicomplex
may be viewed as a completed cobar construction.

\medskip\noindent{\bf
3. Representations of free simplicial groupoids}\smallskip\noindent
Recall
that any cosimplicial manifold
$M=\{M_{\sharp}\}$ gives rise to a simplicial
differential graded de Rham algebra
$
\Omega M = (\Omega^*(M_{\sharp}),d, \dots),
$
cf. \cite\bottsega.
Here $\Omega^j(M_q)$
are the $j$-forms on $M_q$, for $q \geq 0$,
the operator $d$ is the usual de Rham operator on each
$M_q$,
and $\dots$ stands for the operators
between usual de Rham algebras
induced by the cosimplicial structure.
In particular,
let
$\Omega \Cal H$
be the
simplicial differential graded algebra
of de Rham forms on
the cosimplicial manifold
$\Cal H = \roman{Hom}(K,G)$.
The construction in the previous Section
yields de Rham forms
on $\Cal H$.
We shall explain this in Section 4 below.
We need some preparation first.
\smallskip
Let $k \geq 0$.
Since $K$ is assumed free,
the product $\roman{Hom}(K,G)
\times K^k$
inherits a canonical structure
of {\it cosimplicial-simplicial\/}
manifold,
and the {\it coend\/}
\linebreak
$
\roman{Hom}(K,G)\times_{\Delta}K^k $,
cf. e.~g. \cite\maclbotw,
is a (non-connected) smooth manifold;
actually the coend will
play no r\hataccent ole in this paper.
Moreover,
the canonical evaluation map
$$
E \colon
\roman{Hom}(K,G)
\times K^k
@>>> G^k
\tag3.1
$$
is well defined and smooth in the  sense
that, for each (simplicial) degree $q$, the corresponding
component
$$
E_q \colon
\roman{Hom}(K_q,G)
 \times  K_q^k
@>>> G^k
\tag3.2
$$
is smooth;
we note that the evaluation map factors through
the coend
\linebreak
$\roman{Hom}(K,G) \times_{\Delta}K^k$
but this will not be important for us.
We can now apply the construction in the previous Section
separately for each simplicial degree $q$.
However the naturality of the constructions
provides mores structure.
\smallskip
Write
$C_{\natural}(\Pi)$
and
$C^{\natural}(\Pi)$
for the {\it normalized\/} Eilenberg-Mac Lane
chains and cochains, respectively, of an ordinary groupoid
$\Pi$ so that
$C_{\natural}(\Pi) = C_{\natural}(N\Pi)$
and
$C^{\natural}(\Pi) = C^{\natural}(N\Pi)$.
The {\it nerve\/} or {\it simplicial bar construction\/}
$NK$ of $K$ inherits a structure of bisimplicial set,
one simplicial structure coming from
that of $K$ and the other one from the nerve construction.
Its bicomplex $CNK$ of
chains
which are {\it normalized\/}
in the $\natural$-direction
looks like
$$
(C_\natural(K_\sharp), \partial_\natural,\partial_\sharp).
\tag3.3
$$
Its
vertical differentials
$\partial_\sharp$
are induced by
the alternating sums of the
face operations
induced by the simplicial structure of $K= \{K_\sharp\}$
and its
horizontal ones
$\partial_\natural$
by
the alternating sums of the
face operations
induced by the nerve construction for $K_q$
separately for each $K_q$.
Normalizing
in the $\sharp$-direction
yields the bicomplex
$$
(\overline C_\natural(K_\sharp), \partial_\natural,\partial_\sharp)
\tag3.4
$$
where
the notation
$\partial_\natural,\partial_\sharp$
is abused.
Its total complex
$$
|NK| =(|\overline CNK|,\partial) =
|(\overline C_\natural(K_\sharp), \partial_\natural,\partial_\sharp)|
$$
has
$|NK|_0 = \bold Z$ and
$$
|NK|_r = \overline C_r(K_0) \oplus \overline C_{r-1}(K_1) \oplus \dots
\oplus \overline C_1(K_{r-1}),\quad r \geq 1,
\tag3.5
$$
and the total differential
$\partial$
is given
by
$$
\partial = \partial_\natural +
 \partial_\flat
\tag3.6
$$
where
on elements of $\overline C_k(K_{r-k})$
the operator $\partial_\flat$ looks like
$$
\partial_\flat
= (-1)^k \partial_\sharp
\colon
\overline C_k(K_{r-k})
@>>>
\overline C_k(K_{r-k-1})  .
\tag3.7
$$
Note that, by normalization,
$\overline C_0(K_r)$ is zero for $r \geq 1$.
Thus for
$$
c = c_{r,0} + c_{r-1,1} + \dots + c_{1,r-1},
\quad
c_{k,q} \in \overline C_k(K_q),
\tag3.8
$$
we have
$$
\aligned
\partial(c) &= \partial_\natural(c) +
\sum_{k+q=r} (-1)^k \partial_\sharp c_{k,q}
\\
&=
\quad\quad
\partial_\natural(c_{r,0})
+
(-1)^{r-1} \partial_\sharp(c_{r-1,1})
\\
&
\quad
+ \partial_\natural(c_{r-1,1})
+
(-1)^{r-2} \partial_\sharp(c_{r-2,2})
\\
&
\quad
+
\cdots
\\
&
\quad
+
\partial_\natural(c_{2,r-2})
+ \partial_\sharp(c_{1,r-1}) .
\endaligned
\tag3.9
$$
\smallskip
Let $K=KY$ for a CW-complex $Y$.
Since $K$ is a loop complex
for $Y$,
the map from $Y$ to
$B|K|$ is a homotopy equivalence,
and the
two spaces
$B|K|$
and
$|NK|$
are homeomorphic;
as CW-complexes they are not the same, though,
and the cell decomposition
of
$|NK|$
must be refined,
in the same way as the canonical homeomorphism
between
the realization $|S_1 \times S_2|$
of the
product of two simplicial sets $S_1$ and
$S_2$
and
the product
$|S_1| \times |S_2|$
of the realizations
will be a cellular
isomorphism only after refinement of the decomposition
of
$|S_1 \times S_2|$;
cf. \cite\puppeone\
for details.
It follows
that the homology of $|NK|$ coincides with that of $Y$.
However this may be seen directly.
To this end we observe at first that,
for  $\sharp$
fixed,
since each $K_\sharp$ is a free group,
the chain
complex (3.3)
amounts
to an exact sequence
$$
0
@<<<
\roman H_1(K_\sharp)
@<{\varepsilon_\sharp}<<
C_1(K_\sharp)
@<<<
\cdots
@<<<
C_k(K_\sharp)
@<<<
\cdots
\tag3.10
$$
and (3.3),
viewed as a simplicial chain complex,
with simplicial structure
in the $\sharp$-direction,
induces
a structure of simplicial abelian group
$
\roman H_1(K_\sharp)= \{\roman H_1(K_q)\}_{q \geq 0}.
$
Here we have written
$\varepsilon_\sharp$
for the projection from
1-cycles to homology.
For $q \geq 0$, denote by $\overline K_q$
the
free group generated by the
degree $q$ generators,
that is, by the non-degenerate
basis elements of $K_q$.
By construction,
the normalized chain complex
$
|\roman H_1(K_\sharp)|
$
of
$\roman H_1(K_\sharp)$
has
$$
|\roman H_1(K_\sharp)|_q
=
\roman H_1(\overline K_q)
=
(\overline K_q)^{\roman {Ab}}
=
C_{q+1}(Y),
\quad
q \geq 0,
\tag3.11
$$
where
$C_*(Y)$ refers to the
cellular chains
on $Y$.
Furthermore,
$
\roman H_0(K_\sharp)
=\{\roman H_0(K_q)\}_{q \geq 0}
$
amounts to the free simplicial abelian  group
generated by a single point.

\proclaim{Proposition 3.12}
The canonical projection map from
{\rm (3.3)}
onto
$\roman H_1(K_\sharp)$
induced by $\varepsilon_\sharp$
together with the canonical map from
$|NK|_0 = \bold Z$
onto
$C_0(Y) = \bold Z$
passes to a deformation retraction
from
$|NK|$
onto $C_*(Y)$
which is natural in $Y$.
\endproclaim

\demo{Proof}
The canonical projection map from
{\rm (3.3)}
onto
$\roman H_1(K_\sharp)$
induced by $\varepsilon_\sharp$
yields a deformation retraction
from  the totalization
$|(C_\natural(K_\sharp), \partial_\natural,\partial_\sharp)|$
onto
the totalization
of
$\roman H_1(K_\sharp)$.
A little thought
reveals that this implies the claim. \qed \enddemo

\smallskip\noindent
N. B.
The
normalization
$\overline C_\natural(K_\sharp)$
contains
$C_\natural(\overline K_\sharp)$
in an obvious fashion but
does {\it not\/}
coincide with
it since products of degenerate free generators
are in general non-degenerate.

\medskip\noindent{\bf 4. Forms on representations of free simplicial groupoids}
\smallskip\noindent
We can now extend
the construction of forms in Section 2
to representations of the simplicial groupoid $K$.
For each simplicial degree $q$,
with $\Pi = K_q$, the pairing (2.12)
looks like
$$
\left(\Omega_G^{*,*}(G^{\natural}); d,\delta_G,\delta^{\natural}\right)
\otimes (C_\natural(K_q),\partial_\natural)
@>>>
\left(\Omega_G^{*,*}(\roman H_q); d,\delta_G\right)
\otimes (C_\natural(K_q),\partial_\natural),
\tag4.1
$$
and these assemble to the pairing
$$
\left(\Omega_G^{*,*}(G^{\natural}); d,\delta_G,\delta^{\natural}\right)
\otimes (C_\natural(K_\sharp),\partial_\natural)
@>>>
\left(\Omega_G^{*,*}(\roman H_\sharp); d,\delta_G\right)
\otimes (C_\natural(K_\sharp),\partial_\natural).
\tag4.2
$$
The left- and right-hand side of (4.2)
both
inherit a simplicial structure
from that of
$K$;
in fact,
on the left-hand side we have such a structure
on $C_\natural (K)$ and, on the right-hand side,
the induced cosimplicial structure
on $\Cal H = \roman{Hom}(K,G)$
induces a simplicial structure
on
$\left(\Omega_G^{*,*}(\Cal H); d,\delta_G\right)$.
The naturality of the constructions implies that (4.2) is compatible with
these
structures, whence we arrive at the pairing
$$
\left(\Omega_G^{*,*}(G^{\natural}); d,\delta_G,\delta^{\natural}\right)
\otimes (C_\natural(K_\sharp);\partial_\natural,\partial_\sharp)
@>>>
\left([(\Omega_G^{*,*}(\roman H_\sharp); d,\delta_G)
\otimes (C_\natural(K_\sharp),\partial_\natural)], \partial_\sharp\right)
\tag4.3
$$
compatible with all the operators.
In a given quintuple degree $(i,j,k,k+\ell,q)$, this pairing
goes from
$
\Omega_G^{i,j}(G^{k})
\otimes C_{k+\ell}(K_q)
$
to
$
\Omega_G^{i,j}(\roman H_q)
\otimes C_{\ell}(K_q)$.
We note that, with reference to the $\sharp$-grading,
$
(\Omega_G^{*,*}(\roman H_\sharp); d,\delta_G)
\otimes (C_\natural(K_\sharp),\partial_\natural)
$
is the
graded object underlying the
diagonal of a certain bisimplicial object;
we have chosen the parentheses \lq $[$\rq\ and
\lq $]$\rq\
on the right-hand side of (4.3),
with the operator
$\partial_\sharp$ outside
these parentheses
to indicate this.
\smallskip
The {\it normalization\/}
${
(\overline\Omega_G^{*,*}(\roman H_\sharp); d,\delta_G,\partial_\sharp)
=
\left(\Omega_G^{*,*}(\roman H_\sharp)
\big /\Omega_G^{*,*}(\roman H_\sharp)^{\roman{degen}};
 d,\delta_G,\partial_\sharp\right)
}$
of
\linebreak
$(\Omega_G^{*,*}(\roman H_\sharp); d,\delta_G,\partial_\sharp)$
is the quotient by the
subspace
of degenerates
where,
for each simplicial degree $q \geq 1$,
the subspace
$
\Omega_G^{*,*}(\roman H_q)^{\roman{degen}}
$
of degenerates
is the sum of the images of the degeneracy operations
$
s_j
$
from
$\Omega_G^{*,*}(\roman H_{q-1})
$
to
$\Omega_G^{*,*}(\roman H_q)$,
for
$0 \leq j \leq q-1$.
Here the notation $d,\delta_G,\partial_\sharp$ is abused.
Ignoring the equivariant theory for the moment, we
recall \cite\bottsega\
that the {\it realization\/}
$
(|\Omega(\Cal H)|,D)
$
of
$\Omega(\Cal H) =(\Omega^*(\roman H_\sharp),d,\partial_\sharp)$
is the total cochain complex of the normalized bicomplex
$$
\Omega^*(\roman H_{0})
@<{\partial_\sharp}<<
\overline{\Omega^*(\roman H_1)}
@<{\partial_\sharp}<<
\dots
@<{\partial_\sharp}<<
\overline{\Omega^*(\roman H_q)}
@<{\partial_\sharp}<<
\dots
$$
whose
vertical differentials are the de Rham operators and whose
horizontal ones
$\partial_\sharp$
are
induced by
the alternating sums of the simplicial
operations
$\partial_p\colon
\Omega^*(\roman H_{q})
@>>>
\Omega^*(\roman H_{q-1})$.
The
graded module
$|\Omega(\Cal H)|$
underlying the {\it total
complex\/}
$$
|(\Omega^*(\Cal H),d,\partial_\sharp )|
=(|\Omega(\Cal H)|,D)
$$
of this bicomplex is by definition
in degree $r$ the direct {\it sum\/}
(not the product)
of the
$\overline{\Omega^p(\roman H_{q})}$
for $p-q=r$.
Thus
$$
|\Omega(\Cal H)|^r
=
\Omega^r(\roman H_0)
\oplus
\overline{\Omega^{r+1}(\roman H_{1})}
\oplus
\cdots
\oplus
\overline{\Omega^{r+q}(\roman H_{q})}
\oplus
\dots
$$
Notice when
$K$ has a finite set of free generators this sum is finite in each degree.
Moreover this construction has an obvious extension
$$
|\Omega_G(\Cal H)|=
(|\Omega_G^{*,*}(\Cal H)|,D)
=
\left|(\Omega_G^{*,*}(\roman H_\sharp);d,\delta_G,\partial_\sharp)\right|
$$
to the equivariant theory
so that
$$
|\Omega_G^{*,*}(\Cal H)|^r
=
|\Omega_G^{*,*}|^r(\roman H_0)
\oplus
\overline{|\Omega_G^{*,*}|^{r+1}(\roman H_{1})}
\oplus
\cdots
\oplus
\overline{|\Omega_G^{*,*}|^{r+q}(\roman H_{q})}
\oplus
\dots
$$
\smallskip
The compatibility of (4.3) with all the operators entails that
after totalization
and normalization we
arrive at the pairing
$$
\left|
\left(\Omega_G^{*,*}(G^{\natural}); d,\delta_G,\delta^{\natural}\right)
\right|
\otimes
\left|
NK
\right|
@>>>
\roman{N}_\sharp
\left(
[
|(\Omega_G^{*,*}(\roman H_\sharp); d,\delta_G)|
\otimes
(C_\natural(K_\sharp),\partial_\natural)
],\partial_\sharp\right)
\tag4.4
$$
where $\roman{N}_\sharp$
refers to normalization in the $\sharp$-direction.
The
generalized
Eilenberg-Zilber theorem
\cite\doldpupp\
yields a natural chain equivalence
from
the right-hand side
of (4.4)
onto
$(|\Omega_G^{*,*}(\Cal H)|,D)
\otimes
|NK|
$.
Hence (4.4) combined with this surjection yields
the pairing
$$
\left|\left(\Omega_G^{*,*}(G^{\natural});
d,\delta_G,\delta^{\natural}\right)
\right|
\otimes
|NK|
@>>>
(|\Omega_G^{*,*}(\Cal H)|,D)
\otimes
|NK|
\tag4.5
$$
When
we combine it with the chain map
$\roman{Id} \otimes \varepsilon$
where
$\varepsilon$
is the augmentation map
from
$|NK|$
to the reals
induced by the obvious projection from $NK$ to a point,
we arrive at the pairing
$$
\langle\cdot,\cdot\rangle \colon
\left|\left(\Omega_G^{*,*}(G^{\natural});
d,\delta_G,\delta^{\natural}\right)
\right|
\otimes
|NK|
@>>>
(|\Omega_G^{*,*}(\Cal H)|,D)
\tag4.6
$$
This pairing can be understood without explicit reference
to the generalized Eilenberg-Zilber theorem:
it amounts to picking the components of the right-hand side
of (4.4) which involve only
$C_0(K_\sharp)$
and ignoring the rest, but the generalized Eilenberg-Zilber
theorem provides the appropriate formal circumstances.
The precise geometric analogue
$$
\Omega_G^{*,*}(BG) \otimes C_*(M)
@>>>
\Omega_G^{*,*}(\roman{Smooth}^o(M,BG) )
$$
of (4.6)
for a smooth manifold $M$
arises from the evaluation pairing
from
\linebreak
$\roman{Smooth}^o(M,BG) \times M$
 to $BG$
combined with integration against chains on $M$ and
subsequent composition with the chain map induced
by the augmentation map from $C_*(M)$ to the reals $\bold R$.
\smallskip
The pairing (4.6) produces
equivariant forms
on the
realization
$\left|(\overline\Omega_G^{*,*}(\roman H_\sharp); d,\delta_G,\partial_\sharp)
\right|$
and hence, as we shall see later,
on the cosimplicial manifold $\Cal H =\roman{Hom}(K,G)$,
in the following way:
Let $u$ and $r$ be
positive integers, let
$
\Omega
$ be an equivariant
form, that is, an element of
$$
\left|(\Omega_G^{*,*}(G^\natural); d,\delta_G,\delta^\natural)\right|
=
\left(|\Omega_G^{*,*}(G^\natural)|,D\right)
$$
of total degree
$u+r$,
and
let $c$
be a  chain
of $|NK|$
of total degree $r\geq 1$,
cf. (3.5).
For $k$ fixed,
let $\Omega^k \in \oplus_{i+j =u+r-k} \Omega_G^{i,j}(G^k)$
be the indicated component.
Then
$$
\langle
\Omega, c
\rangle=
\langle
\Omega^1, c_{1,r-1}
\rangle
+
\dots
+
\langle
\Omega^{r}, c_{r,0}
\rangle
\in
|\overline \Omega_G^{*,*}(\roman H_\sharp)|.
$$

\proclaim{Lemma 4.7}
Suppose $\Omega$ and $c$ are closed. Then
$\langle
\Omega, c
\rangle$
is a closed form in
$$
\left(\left|\overline \Omega_G^{*,*}(\roman H_\sharp)\right|,D\right)
=
\left|\left(\overline \Omega_G^{*,*}(\roman H_\sharp),
d,\delta_G,\partial_\sharp \right)\right|
$$
of total degree $u$.
\endproclaim

\demo{Proof}
This follows at once from the identity
$$
D\langle \Omega,c\rangle
=
\langle d_G \Omega,c\rangle
+ (-1)^{|\Omega|}
\langle \Omega,\partial c\rangle .
$$
\enddemo

\beginsection 5. Realization and integration

A standard construction endows
$
(|\Omega^{*}(\Cal H)|,D)
$
with a structure of differential graded algebra;
we shall explain this below in the equivariant setting.
Before doing so
for illustration,
we point out that, for $K$ the
Kan group on the 2-sphere with its standard decomposition
with two cells,
for a Lie group $G$,
$|\roman{Hom}(K,G)|$
amounts to the space $\Omega G =\roman{Map}^o(S^1,G)$ of
based
loops on $G$
and we have on the one hand the bar construction
$B \Omega^*G$ on the de Rham complex
$\Omega^*G$,
with its shuffle multiplication,
 as a model for the algebra of cochains on
the based loop space.
On the other hand,
as a cosimplicial space,
$\Cal H =\roman{Hom}(K,G)$
looks like
$(o,G,G^2, \dots, G^q,\dots)$,
and the realization of the simplicial differential graded
algebra
$\Omega^*(\Cal H)$
(whose algebra structure is yet to be explained)
has components
$\Omega^*(G^q)$.
The canonical maps
from
$(\Omega^*(G))^q$
to $\Omega^*(G^q)$
now induce a morphism of differential graded algebras
from
$B \Omega^*G$
to
$
(|\Omega^*(\Cal H)|,D)
$
which is a homology isomorphism.
\smallskip
We now recall the
construction
of differential graded algebra
structure on
$
(|\Omega_G^{*,*}(\Cal H)|,D)
$:
For each pairs $(i,j)$ and
$(i',j')$
of bidegrees,
we have the simplicial vector spaces
$\Omega_G^{i,j}( \roman H_\sharp)$
and
$\Omega_G^{i',j'}( \roman H_\sharp)$,
and for each pair $(q,q')$ the {\it shuffle map\/}
$\nabla$ yields
a natural morphism
$$
\nabla
\colon
\Omega_G^{i,j}( \roman H_q)
\otimes
\Omega_G^{i',j'}( \roman H_{q'})
@>>>
\Omega_G^{i,j}( \roman H_{q+q'})
\otimes
\Omega_G^{i',j'}( \roman H_{q+q'})
$$
of vector spaces which, combined with usual multiplication
of forms, yields a pairing
$$
\Omega_G^{i,j}( \roman H_q)
\otimes
\Omega_G^{i',j'}( \roman H_{q'})
@>>>
\Omega_G^{i+i',j+j'}( \roman H_{q+q'}).
$$
This
pairing endows
$
(|\Omega_G^{*,*}(\Cal H)|,D)
$
with a structure of differential
graded commutative algebra
which is natural in the data;
by construction, it arises from a differential
trigraded algebra structure.
\smallskip
We now relate this algebra
with forms on the geometric realization.
To this end,
pick $q \geq 0$ and
consider the evaluation mapping
from
$
\Delta_q \times \roman{Smooth}(\Delta_q, \roman H_q)
$
to
$\roman H_q$.
For each $r \geq 0$,
integration over $\Delta_q$
induces a map
$$
I_q
\colon
\Omega^{r+q}(\roman H_{q})
@>>>
\Omega^r(\roman{Smooth}(\Delta_q, \roman H_q))
$$
with a suitable interpretation of forms
on the mapping spaces.
Using the theory of {\it differentiable space\/}
\cite\chen, \cite\chentwo\ or, what amounts to the same,
that of {\it \lq\lq diffeological\rq\rq\  space\/}
(\lq\lq espace diff\'eologique\rq\rq)
\cite\sourithr, \cite\sourifiv,
forms
on the mapping spaces
admit a purely
finite dimensional interpretation
and do {\it not\/} require infinite dimensional techniques.
With a suitable interpretation of forms
$\Omega^*(|\Cal H|_{\roman{smooth}})$ on the geometric realization
$|\Cal H|_{\roman{smooth}}$,
the
integration maps assemble to  a morphism
$$
I
\colon
|(\Omega^*(\Cal H); d,\partial_\sharp)|
=
(|\Omega^*(\Cal H)|, D)
@>>>
(\Omega^*(|\Cal H|_{\roman{smooth}}),d)
\tag5.1
$$
of differential graded algebras,
cf. Section 5 of \cite\bottsega.
The
problem here is that
the geometric realization $|\Cal H|_{\roman{smooth}}$
will have singularities.
This difficulty is overcome by means of the already cited
concept of differentiable space, in the following way:
Recall that a {\it plot\/}
for
$|\Cal H|_{\roman{smooth}}$
is a map $F$ from a smooth finite dimensional manifold
$W$ to
$|\Cal H|_{\roman{smooth}}$
which is smooth in the sense that the adjoint
$$
F_q^\natural
\colon W\times \Delta_q @>>> \roman H_q
$$
of each component
$$
F_q
\colon W  @>>> \roman{Smooth}(\Delta_q,\roman H_q)
$$
of $F$ is smooth
\cite{\chen, \chentwo};
in this theory,
a {\it form\/}
on
$|\Cal H|_{\roman{smooth}}$
is the assignment of a form on $W$ to each plot
which is natural for smooth maps in the domains of the plots.
With this interpretation
of forms
on
$(\Omega^*(|\Cal H|_{\roman{smooth}}),d)$
the above integration mapping $I$ makes strict sense.
In particular, when $\Cal M$ is a subspace of
$|\Cal H|_{\roman{smooth}}$
which {\it is\/} smooth
we can combine the integration mapping with restriction,
and there results a morphism
$I$
of differential graded algebras
from
$(|\Omega^*(\Cal H)|, D)$ to
$(\Omega^*(\Cal M),d)$.
Furthermore the whole construction is $G$-invariant whence,
with the appropriate notion of $G$-equivariant plots,
we finally obtain a morphism
$$
I
\colon
|(\Omega_G^{*,*}(\Cal H); d,\delta_G,\partial_\sharp)|
=
(|\Omega_G^{*,*}(\Cal H)|, D)
@>>>
|(\Omega_G^{*,*}(|\Cal H|_{\roman{smooth}}); d,\delta_G)|
\tag5.2
$$
of differential bigraded algebras.
\smallskip
Under our circumstances,
$G$-equivariant plots admit a natural interpretation
as {\it $G$-equivariant families of principal bundles with connection\/}:
a $G$-equivariant plot
$F \colon W \to |\Cal H|_{\roman{smooth}}$,
combined with the map $\Phi$
(cf. (1.6)),
yields a
map from $W$ to $\roman{Smooth}^o(Y,BG)$
having a \lq\lq smooth\rq\rq\
$G$-equivariant
adjoint
$$
\widetilde F \colon W \times Y
@>>>
BG
$$
satisfying
$\widetilde F(w,o) = o$.
Consequently
a $G$-equivariant plot
$F$ for
$|\Cal H|_{\roman{smooth}}$
defined on $W$
amounts to a smooth $G$-equivariant family
of $G$-bundles with connection on $Y$
parametrized by $W$.
\smallskip
A cosimplicial space
 is said to {\it converge\/}
\cite\anderone, \cite\bottsega,
\cite\bouskan,
when integration yields a cohomology equivalence
from the cohomology of the realization
of the simplicial de Rham algebra
to the cohomology of the realization.
The cosimplicial space $\Cal H$ will rarely have this property;
however see Section 7 below.

\medskip\noindent{\bf 6. Extended moduli spaces for a closed surface and
generalizations}\smallskip\noindent
Let $\Sigma$ be a closed topological surface
of genus $\ell \geq 0$,
endowed with the usual CW-decomposition
with a single 0-cell $o$,
with 1-cells
$u_1,v_1,\dots,u_\ell,v_\ell$, and with a single 2-cell
$c$,
and let
$$
\Cal P = \langle
x_1,y_1,\dots,x_\ell,y_\ell; r\rangle
$$
be the corresponding presentation for the fundamental group $\pi$ of $\Sigma$.
We maintain the notation as in Section
2 of \cite\huebkan\
without repeating it.
Our aim is
to show how the results of
\cite\modus, \cite\huebjeff\
may at once be deduced
from our general theory:
Write
$Q$ for the given bilinear 2-form on $g$,
and
let
$$
\Omega_Q =
Q^{0,3,1}+ Q^{2,1,1}
+
Q^{0,2,2}
+
Q^{2,0,2} \in \left|\Omega_G^{*,*}(G^{\natural})\right|
$$
be the resulting equivariantly closed form
of total degree 4,
cf. (2.17).
Since $\roman H_2(\Sigma)$ is infinite cyclic,
in view of (3.12), there is a 2-cycle
$$
c = c_{2,0} + c_{1,1}, \quad
c_{2,0} \in C_2(K_0),
\ c_{1,1} \in \overline C_1(K_1),
$$
of $|NK|$
which under the deformation retraction
onto the cellular chains of $\Sigma$
goes to a 2-cycle representing a generator.
It is very easy to manufacture such a 2-cycle:
Let
$c_{2,0} \in C_2(K_0)= C_2(F)$
be a 2-chain with
$\partial_\natural c_{2,0} =
\Pi [x_j,y_j] \in F$;
such a
$c_{2,0}$ exists since
$\Pi [x_j,y_j]$ is zero in $\roman H_1(F) = F^{\roman{Ab}}$;
moreover,
let
$c_{1,1} = r \in K_1$ so that,
by construction,
$\partial_\sharp  c_{1,1} =
\Pi [x_j,y_j] \in K_0$. Then
$c$ is closed in $|NK|$, and
$
\langle Q,c\rangle
$
is a closed element of
$\left|(\Omega_G^{*,*}(\roman H_\sharp); d, \delta_G, \partial_\sharp)\right|$
of degree 2.
Notice when $\ell =0$ we have
$c_{2,0} = 0$.
\smallskip
Embed the Lie algebra $g$ into
$\roman{Smooth}(\Delta_1,G)$
by the assignment to $X \in g$ of
the corresponding path $t \mapsto \roman{exp}(tX)$,
let $O \subseteq g$ be the subspace where the exponential
mapping is regular, and let
$\Cal M$ be the subspace
of $|\Cal H|_{\roman{smooth}}$
consisting of pairs $(w,X) \in G^{2\ell} \times O $ so that
$\roman{exp}(X) = r(w)$.
This is a smooth finite dimensional
$G$-manifold
and the inclusion $F$ from $\Cal M$
to
$|\Cal H|_{\roman{smooth}}$
is a $G$-equivariant plot.
By construction, the equivariantly closed form
$$
I\langle Q,c\rangle
\in \left|\overline \Omega_G^{*,*}(\Cal M)\right|
$$
of degree 2 has components
$$
\align
\omega_c &=
I\langle Q^{0,3,1},c_{1,1}\rangle
+
I\langle Q^{0,2,2},c_{2,0}\rangle
 \in \Omega_G^{0,2}(\Cal M)
\\
\mu^{\sharp}&=
I\langle Q^{2,1,1},c_{1,1}\rangle
+
I\langle Q^{2,0,2},c_{2,0}\rangle
\in \Omega_G^{2,0}(\Cal M)
\endalign
$$
so that
$
\mu^{\sharp}
\colon g \to C^{\infty}(\Cal M)
$
is the adjoint of a smooth map
$\mu$ from $\Cal M$ to $g$.
The
sum
$\omega_c + \mu^{\sharp}$
to be equivariantly closed
amounts
to the closedness of $\omega_c$ in the usual sense
together with the property that
$$
\delta_G \omega_c = d \mu^{\sharp}
$$
which
is the momentum mapping property.
In particular,
the integration mapping
from the realization of the equivariant
simplicial differential graded de Rham algebra
to the equivariant differential graded
de Rham algebra on the realization
reproduced in Section 5 above
now provides
a natural explanation
for the operation of integration
along linear paths in $g$
which in \cite\modus, \cite\huebjeff\
seemed somewhat ad hoc.
The term
$I\langle Q^{2,0,2},c_{2,0}\rangle$
is actually irrelevant and may be ignored;
it amounts to a constant modification of the momentum mapping.
\smallskip
The whole approach may be extended to
arbitrary 2-complexes $Y$ with a single $0$-cell and
a single 2-cell;
write $r$ for the corresponding relator.
The fundamental group $\pi$ of such a 2-complex
is a one-relator group.
It is well known that $Y$ and $\pi$ have
second homology group a copy of the integers
if and only if the exponent sum of each generator
in the relator $r$
equals zero.
In this case the above construction carries over verbatim
and there results an extended moduli space
$\Cal M$ together with a closed 2-form
$\omega$ and momentum mapping.
However in order for $\omega$
to be non-degenerate
we must require the relevant cup pairings
on $\roman H^1(\pi,\cdot)$
to be non-degenerate.
This is related with the question whether
$\pi$ is a 2-dimensional Poincar\'e duality
group over the reals.
We do not pursue this question here.

\medskip\noindent{\bf 7. Cohomology}\smallskip\noindent
Let $r \geq 1$, let $c$ be a cellular $r$-cycle of $Y$ representing
an integral homology class, and let
$c_K$ be an $r$-cycle
of $|NK|$
which under the deformation retraction
onto the cellular chains of $Y$
goes to $c$, cf. (3.12).
For every
invariant polynomial $Q$ on $g$ of degree $u$,
with $r \leq 2u$,
by (4.7),
$\langle \Omega_Q,c_K\rangle$
is a closed element of
$
(|\Omega_G^{*,*}(\Cal H)|,D)
$
of degree $2u-r$;
here $\Omega_Q$
refers to the corresponding closed elements of
$\left|(\Omega_G^{*,*}(G^{\natural});d,\delta_G,\delta^{\natural})\right|$,
cf. (2.17).
Recall that in Section 1 we constructed a weak chain equivalence
between the space of based maps from $Y$ to $BG$ and
the realization of $\Cal H$.

\proclaim {Theorem 7.1}
As a graded commutative algebra,
the equivariant cohomology
of each connected component of
$|\Cal H|$ and hence of
the space of based
maps from $Y$
to $BG$
is freely generated by
the classes of
the elements $I\langle \Omega_Q,c_K\rangle$
where $Q$ runs through a set of invariant polynomials
on $g$ and $c$
through a set of representatives in degree $\geq 1$
of the real homology of $Y$ subject to the restriction
$|c|< 2 |Q|$,
together with the invariant polynomials $Q$ viewed
as elements of
$\Omega_G^{*,0}(G^0)$.
\endproclaim

\demo{Proof}
This is proved by induction on dimension,
with reference to the fibration
(1.8.1).
The argument is formally the same as that
hinted at on p. 181 of \cite\donalkro, cf. also
Note 5.1.2 on p. 206.
The induction starts with the observation
that
$\roman{Hom}(KY^1,G)$
amounts to a product of as many copies of $G$ as $Y$ has 1-cells
and that, for a circle $K=S^1$,
when $c$ represents the generator of its first
homology group,
the element
$\langle \Omega_Q,c_K\rangle$
yields
the exterior generator of
$\roman H^{2|Q|-1}(G)$
which transgresses
to the class of $[Q]$ in
$\roman H^{2|Q|}(BG)$.
We leave the details to the reader. \qed
\enddemo

Here
the realization
$|\Cal H|_{\roman{smooth}}$
of $\Cal H$
is viewed as a  space
with differentiable
structure
in the sense
of \cite{\chen, \chentwo, \sourithr, \sourifiv}
as explained
in Section 5 above.
\smallskip\noindent
{\smc Illustration.}
Let $Y$ be a closed surface $\Sigma$ of genus $\ell \geq 0$,
let $G = \roman U(n)$, the unitary group,
let $Q_1, \dots, Q_n$ be the Chern polynomials,
and let
$\Omega_1,\dots,\Omega_n$
be the corresponding closed elements of
$\left|(\Omega_G^{*,*}(G^{\natural});d,\delta_G,\delta^{\natural})\right|$.
Maintaining the notation in the previous Section,
for $r  = 1, \dots, n$ and $j=1,\dots,\ell$,
we get
the elements
$$
\align
f_r &= \langle \Omega_r, c\rangle, \quad |f_r| = 2r-2,
\\
b^j_r &= \langle \Omega_r, u_j\rangle,  \quad |b^j_r| = 2r-1,
\\
b^{j+\ell}_r &= \langle \Omega_r, v_j\rangle,  \quad |b^{j+\ell}_r| = 2r-1,
\\
a_r &= Q_r,\quad |a_r| = 2r,\quad
\text{viewed as an element of\ }\Omega_G^{2r,0}(G^0).
\endalign
$$
They freely generate the equivariant cohomology
of the space $\roman{Map}^o(\Sigma, BG)$
or, what amounts to the same,
of the union
over all topological types of $G$-bundles
of spaces of based
gauge equivalence classes,
cf. what is said in Section 1.
Notice $f_1$ picks the topological type of connected component.
The generators
$f_r$ and $b^j_r$
coincide with those constructed in
\cite\newstone,
see also
Section 2 of \cite\atibottw\
and \cite\jeffrthr.
Likewise,
for genus $\ell = 0$
and arbitrary connected $G$,
when $Q$
denotes
the given invariant
symmetric bilinear form on $g$
so that $Q^{0,3,1}$
is E. Cartan's
fundamental 3-form on $G$,
cf. what was said at the end of Section 2,
the resulting 2-form on $\Omega G$
restricts to the Kirillov
form
on each connected component
of $\roman{Hom}(S^1,G)$,
when identified with the adjoint orbits
in $g$ generated by some $X$ with $\roman{exp}(X) = e$.
\smallskip
When $G$ is simply connected
the cosimplicial space $\Cal H = \roman {Hom}(K\Sigma,G)$
converges, that is,
the integration
mapping (5.2)
is a cohomology equivalence
from the cohomology of the realization
of the simplicial de Rham algebra
to the cohomology of the realization.

\medskip\noindent{\bf 8. 3-complexes and 3-manifolds}
\smallskip\noindent
Let $Y$ be a
3-complex with a single 3-cell,
for example,
a closed compact 3-manifold,
endowed with a regular CW-decomposition
with a single 0-cell $o$,
with 1-cells
$u_1,\dots,u_\ell$, 2-cells
$c_1,\dots,c_\ell$, and a single
3-cell $c$.
We maintain the notation  in  Section
3 of \cite\huebkan\
without repeating it.
We shall give
a purely finite dimensional description of the Chern-Simons function.
This will be the assignment
to
every $G$-equivariant plot
$$
F
\colon
W
@>>>
|\Cal H|_{\roman{smooth}}
$$
of a smooth
$G$-invariant
map $\Psi$ from $W$
to the circle $S^1$
which is natural in $G$-equivariant plots.
\smallskip
Write
$Q$ for the given bilinear 2-form on $g$,
and
let
$$
\Omega_Q =
Q^{0,3,1}+ Q^{2,1,1}
+
Q^{0,2,2}
+
Q^{2,0,2} \in \left|\Omega_G^{*,*}(G^{\natural})\right|
$$
be the resulting equivariantly closed form
of total degree 4,
cf. (2.17).
Suppose
$\roman H_3(Y)$ infinite cyclic;
for example this will be true when
$Y$ is an orientable 3-manifold.
In view of (3.12), there is a 3-cycle
$$
c = c_{3,0} +c_{2,1} + c_{1,2}, \quad
c_{3,0} \in C_3(K_0),
\ c_{2,1} \in \overline C_2(K_1),
\ c_{1,2} \in \overline C_1(K_2),
$$
of $|NK|$
which under the deformation retraction
onto the cellular chains of $Y$
goes to a cellular 3-cycle representing a generator
of $\roman H_3(Y)$.
An explicit such a 3-cycle
is obtained as follows:
Let
$c_{1,2} = \sigma\in C_1(K_2)$;
then
$$
\partial_{\sharp}
\sigma
=
(s_0z_1) r^{\varepsilon_1}_{j_1}(s_0z_1)^{-1} \dots
(s_0z_m) r^{\varepsilon_m}_{j_m}(s_0z_m)^{-1} \in C_1(K_1)
$$
and the class of the latter
in
$\roman H_1(K_1)$
is zero.
In fact,
$\roman H_1(K_1)$
is the free abelian group
on the relators
$r_1,\dots,r_\ell$
and the degeneracies
$s_0x_1,\dots,s_0x_\ell$
of the generators
in $\Cal P$,
and the
subgroup generated
by the relators
amounts to
the
group $C_2(Y)$
of cellular 2-chains of $Y$.
The assignment to
$\sigma$
of
the image
of $\partial_{\sharp}
\sigma$
in
$C_2(Y)$
under the
map from
$C_1(K_1)$ to
$\roman H_1(K_1)$
combined with
the projection
onto
$C_2(Y)$,
cf. (3.10) above,
is the value
of the boundary
$\partial \sigma \in C_2(Y)$
under the cellular boundary
$\partial \colon
C_3(Y)
\to
C_2(Y)
$
and this is zero since
$\sigma$
represents a 3-cycle.
However,
the image of $\partial_{\sharp}
\sigma$
in
$\roman H_1(K_1)$
lies in
the subgroup of
$\roman H_1(K_1)$
generated by
the relators
$r_1,\dots,r_\ell$
and this subgroup is mapped isomorphically
onto $C_2(Y)$
whence
the image of $\partial_{\sharp}
\sigma$
in
$\roman H_1(K_1)$
is zero.
Since the sequence (3.10)
is exact, we conclude that there is
a chain
$c_{2,1} \in C_2(K_1)$
with
$$
\partial_{\flat}c_{2,1} =
\partial_{\sharp}c_{1,2} \in C_1(K_1).
$$
Next,
$$
\partial_{\flat}\partial_{\sharp}c_{2,1}
=
\partial_{\sharp}\partial_{\flat} c_{2,1}
=
\partial_{\sharp}\partial_{\sharp}c_{1,2} = 0
$$
whence, again in view of the exactness of (3.10),
there is
a chain
$c_{3,0} \in C_3(K_0)$
with
$$
\partial_{\flat}c_{3,0} =
-\partial_{\sharp}c_{2,1} \in C_2(K_0).
$$
Then
$$
c = c_{3,0} +c_{2,1} + c_{1,2}
$$
is the desired 3-cycle in $|NK|$, and
$
\langle Q,c\rangle
$
is a closed element of
$\left|(\Omega_G^{*,*}(\roman H_\sharp); d, \delta_G, \partial_\sharp)\right|$
of degree 1.
In some more detail,
$
\langle Q,c\rangle
$
has components
$$
\align
\langle Q^{0,3,1},c_{1,2}\rangle
&\in \Omega^{0,3}(\roman H_2),
\\
\langle Q^{0,2,2},c_{2,1}\rangle
 &\in \Omega_G^{0,2}(\roman H_1),
\\
\langle Q^{2,1,1},c_{1,2}\rangle
&\in \Omega^{2,1}(\roman H_2),
\\
\langle Q^{2,0,2},c_{2,1}\rangle
&\in \Omega^{2,0}(\roman H_1).
\endalign
$$
Notice that here
$\langle Q^{0,3,1},c_{1,2}\rangle
\in \Omega^{0,3}(\roman H_2)$
is just the
form pulled back from Cartan's form
$\lambda \in \Omega^3(G)$
via the canonical projection
from $\roman H_2$ onto its primitive part
$P_2 = G$.
Thus,
keeping in mind that
$\roman H_1 =
G^{2\ell}$,
under the present circumstances,
the construction yields
the 2-form
$$
\alpha =\langle Q^{0,2,2},c_{2,1}\rangle
\in \Omega_G^{0,2}(G^{2\ell})
$$
having the property that
$$
d \alpha
= i^* \lambda
\in \Omega_G^{0,3}(G^{2\ell})
$$
where $i$ refers to the smooth map
from
$G^{2\ell}$ to $G$ induced by the
identity among relations
(3.1) in \cite\huebkan.
Notice also that there is no component
involving a form on $\roman H_0$.
This relies on the fact that the (non-equivariant)
Shulman
construction
(2.16)
yields only non-zero forms
in $\Omega^j(G^k)$ for $j \geq k$.
Thus  $(\alpha, \lambda) \in \Omega^2(\roman H_1) \times \Omega^3(G)$
is a pair of forms
which yields
an equivariant form in
$$
|(\Omega_G^{*,*}(\Cal H); d,\delta_G,\partial_\sharp)|
=
(|\Omega_G^{*,*}(\Cal H)|, D)
$$
of total degree 1, and,
under the integration mapping (5.2), this form
passes to an equivariant
1-from
in
$|(\Omega_G^{*,*}(|\Cal H|_{\roman{smooth}}); d,\delta_G)|$.
Hence,
given a $G$-equivariant plot
$$
F
\colon
W
@>>>
|\Cal H|_{\roman{smooth}},
$$
its adjoint $F^{\natural}$ has components
$$
F^{\natural}_q
\colon W \times \Delta_q\to \roman H_q
$$
and here only
$F^{\natural}_1$
and $F^{\natural}_2$
are relevant; they fit into the commutative diagram
$$
\CD
W \times \Delta_1 @>{\roman{Id} \times \varepsilon^2}>>
W \times \Delta_2
\\
@V{F^{\natural}_1}VV
@V{\roman{pr}F^{\natural}_2}VV
\\
\roman H_1
@>>>
G
\endCD
$$
where
$\roman{pr}$
refers to the canonical projection
from $\roman H_2$ onto its primitive part
$P_2 = G$.
Now
$$
\psi =
\int_{\Delta_2} (\roman{pr}F^{\natural}_2)^*\lambda
+
\int_{\Delta_1}(F^{\natural}_1)^*\alpha
$$
is a closed $G$-equivariant 1-form on $W$
having integral periods and hence
integrates to a
smooth
map $\Psi$ from
$W$ to the circle $S^1$.
Moreover,
$\Omega_Q$
is equivariantly
closed,
the term
$Q^{2,0,2}$
is irrelevant,
and, for every $X \in g$,
the value
$\delta_G \psi (X) = -\psi(X_W)$
is calculated by
$$
d I\langle Q^{2,1,1},c_{1,2}\rangle
$$
where
$\langle Q^{2,1,1},c_{1,2}\rangle \in \Omega_G^{2,1}(\roman H_2)$ and
where $d$ is the de Rham operator;
however,
for degree reasons,
$I\langle Q^{2,1,1},c_{1,2}\rangle$
is zero
whence,
for every $X \in g$,
the value
$\delta_G \psi (X) = -\psi(X_W)$
is zero
and hence
$\Psi$ is constant on $G$-orbits,
that is to say, $G$-equivariant.
Thus,
the choice
of the cycle $c$
determines,
for every $G$-equivariant plot
$
F
$ defined on a smooth $G$-manifold $W$,
a smooth
$G$-equivariant
map $\Psi$ from
$W$ to the circle $S^1$
which is natural in $G$-equivariant plots.
This is our description of the Chern-Simons function.
\smallskip
We conclude this Section with an observation:
Formally,
we can interpret the
 smooth map
$i$
from
$G^{2\ell}$ to $G$ induced by the above
mentioned identity among relations
 as arising from the presentation
$$
\widetilde{\Cal P} = \langle
x_1,\dots,x_\ell, r_1,\dots,r_\ell; s\rangle,
\quad
s =
z_1 r^{\varepsilon_1}_{j_1}z_1^{-1} \dots
z_m r^{\varepsilon_m}_{j_m}z_m^{-1}
$$
of a one relator group arising
from $K_1$,
the free group on
$x_1,\dots,x_\ell, r_1,\dots,r_\ell$,
by interpreting the
identity among relations
(3.1) in \cite\huebkan\
as a relation among the generators of
$K_1$.
The map $i$ from
from
$G^{2\ell}$ to
$G$
then amounts to the word map
given by the association
$$
(a_1,\dots,a_\ell, b_1,\dots, b_{\ell})
\longmapsto
w_1(a) b^{\varepsilon_1}_{j_1}(w_1(a))^{-1} \dots
w_m(a) b^{\varepsilon_m}_{j_m}(w_m(a))^{-1}
\in G
$$
where
$w_k(a)\in G$ is obtained by substituting
each occurrence of
$x_j$ in $z_k$ by $a_j$.
Plainly,
then
$$
\align
\omega &= \langle Q^{0,3,1}, c_{1,2}\rangle +
         \langle Q^{0,2,2}, c_{2,1}\rangle
\in \Omega_G^{0,3}(\roman H_2) +\Omega_G^{0,2}(\roman H_1)
\\
\mu^{\sharp} &= \langle Q^{2,1,1}, c_{1,2}\rangle
\in \Omega_G^{2,1}(\roman H_2)
\endalign
$$
satisfy
$$
\delta_G \omega = \pm d \mu^{\sharp}
$$
and $\omega$ is a closed 2-form in the appropriate sense,
just as in the case of surface groups
studied before.

\medskip\noindent{\bf 9. Simply connected 4-manifolds}\smallskip\noindent
As in Section 4 of \cite\huebkan,
we
describe a simply connected 4-manifold $Y$ as
the cofibre of a map $f$
from
the 3-sphere $S^3$ to a bunch $\vee_{\ell} S_j^2$
of $\ell$ copies of the 2-sphere.
We maintain the notation in [ibidem]
but do not reproduce it.
When $Y$ underlies a smooth
4-manifold, the space of based maps
from $Y$ to $BG$
is a model for the space
of based gauge equivalence classes of
connections on all topological types of bundles
on $Y$
and Theorem 7.1 above gives a complete description of its real
equivariant cohomology.
Perhaps
moduli spaces
of
based gauge equivalence classes of
ASD-connections
can be found within
the geometric realization
of our cosimplicial manifold $\roman{Hom}(KY,G)$,
in the following way:
Endow the 4-manifold $Y$ with a metric as usual
and view
$Y$ as a compactification of
the interior $\bold R^4$
of the 4-ball, with the induced metric.
When this metric is flat up to a diffeomorphism
of
$\bold R^4$,
standard constructions
yield all finite energy ASD-connections
on $\bold R^4$.
Whether or not
this happens to be the case,
finite energy ASD-connections
should correspond to certain maps
of the kind $\phi_3$.
In the flat case,
such connections
can be extended
over the 4-sphere,
by Uhlenbeck's
removable singularities theorem.
In general
there are presumably additional constraints
corresponding to the
requirement that
for a choice of $\phi_1$
determined by the data
the diagram
\cite\huebkan\ (4.2)
be commutative.
Perhaps
these additional constraints
explain the additional term
$-db_2^+$
in the formula for the dimension
of the moduli space
which comes from the index theorem
where $d$ refers to the dimension of the structure
group and
$b_2^+$
to the rank of the self dual part of $\roman H^2(Y)$
as usual.
Also in special cases
the description of
ASD-connections
on the 4-sphere
and related 4-manifolds
as
holomorphic maps from the 2-sphere
or related 2-manifolds
to $\Omega G$ might be relevant here,
cf. \cite\atiyahfo,
\cite\akingone.
Another question is
whether anything reasonable can be said
about the component $\phi_1$,
which is a point of the product of $\ell$ loop spaces $\Omega G$.
Is there a way to relate
Yang-Mills theory over $Y$ with Yang-Mills theory
over the embedded 2-spheres
where
critical points
are known to correspond to
geodesics or
homomorphisms?
The map $\Phi$
from the realization of $\Cal H$
to
$\roman{Map}^o(Y,BG)$
\cite\huebkan\
involves choices,
in particular a choice of homotopy inverse
from $Y$ to $|SY|$ and this choice
destroys the symmetry of the situation.
The ASD-condition will presumably only come down to
a certain fixed map from
$S^1$ to $G$ for each 2-sphere in $Y$.
A finer decomposition, e.~g. a triangulation of $Y$,
will
provide a canonical map from
$Y$ to $|SY|$
and hence might restore the missing symmetry.
We hope to settle these issues eventually elsewhere.

\medskip
\widestnumber\key{999}
\centerline{References}
\smallskip\noindent
References \cite\anderone\  -- \cite\wittesix\
are given in the first part \cite\huebkan.
\smallskip\noindent

\ref \no \atiyahfo
\by M. F. Atiyah
\paper Instantons in two and four dimensions
\jour  Comm. Math. Phys.
\vol 93
\yr 1984
\pages 437--451
\endref

\ref \no  \atiboton
\by M. Atiyah and R. Bott
\paper The moment map  and equivariant cohomology
\jour Topology
\vol 23
\yr 1984
\pages  1--28
\endref

\ref \no \bergever
\by N. Berline, E. Getzler, and M. Vergne
\book Heat kernels and Dirac Operators
\publ Springer
\publaddr Berlin $\cdot$ Heidelberg $\cdot$ New York $\cdot$ Tokyo
\yr 1992
\endref

\ref \no \bottone
\by R. Bott
\paper On the Chern-Weil homomorphism and the continuous cohomology of
Lie groups
\jour Advances
\vol 11
\yr 1973
\pages  289--303
\endref

\ref \no \bottsame
\by R. Bott and H. Samelson
\paper Applications of the theory of Morse to symmetric spaces
\jour Amer. J. of Math.
\vol 80
\yr 1958
\pages  954--1029
\endref

\ref \no \botshust
\by R. Bott, H. Shulman, and J. D. Stasheff
\paper On the de Rham theory of certain classifying spaces
\jour Advances
\vol 20
\yr 1976
\pages 43--56
\endref

\ref \no \chen
\by K. T. Chen
\paper Iterated path integrals
\jour Bull. Amer. Math. Soc.
\vol 83
\yr 1977
\pages  831--879
\endref

\ref \no \chentwo
\by K. T. Chen
\paper Degeneracy indices and Chern classes
\jour  Adv. in Math.
\vol 45
\yr 1982
\pages 73--91
\endref
\ref \no \doldpupp
\by A. Dold und D. Puppe
\paper Homologie nicht-additiver Funktoren
\jour Annales de l'Institut Fourier
\vol 11
\yr 1961
\pages   201--313
\endref

\ref \no \guhujewe
\by K. Guruprasad, J. Huebschmann, L. Jeffrey, and A. Weinstein
\paper Group systems,
groupoids, and moduli spaces of parabolic bundles
\paperinfo preprint, 1994
\endref

\ref \no \huebkan
\by J. Huebschmann
\paper
Extended moduli spaces
and the Kan construction
\paperinfo MPI preprint, 1995,
dg-ga/9505005
\endref

\ref \no \akingone
\by A. D. King
\paper Instantons and holomorphic bundles on the blown-up plane
\paperinfo D. Phil. Thesis, Oxford University, April 1989
\endref

\ref \no \newstone
\by P. E. Newstead
\paper Characteristic classes of stable bundles over an algebraic curve
\jour Trans. Amer. Math. Soc.
\vol 169
\yr 1972
\pages 337--345
\endref
\ref \no \sourithr
\by J. M. Souriau
\paper Groupes diff\'erentiels
\paperinfo in: Lecture Notes in Mathematics, No. 836
\publ Springer
\publaddr Berlin-Heidelberg-New York-Tokyo
\yr 1981
\pages 81
\endref

\ref \no \sourifiv
\by J. M. Souriau
\paper Groupes diff\'erentiels
et physique math\'ematique
\paperinfo in: Feuilletages et quantification g\'eom\'etrique,
P. Dazord \& N. Desolneux-Moulis, eds.,
Travaux en cours
\publ Herman
\publaddr Paris
\yr 1984
\endref

\ref \no \witteele
\by E. Witten
\paper Two-dimensional gauge theories revisited
\jour J. Geom. and Physics
\vol 9
\yr 1992
\pages 303--368
\endref

\enddocument